\documentclass[11pt]{article}
\usepackage{latexsym}
\usepackage{amssymb,amsmath}
\usepackage{graphicx}
\RequirePackage{tikz}
\usepackage{epsfig}
\usepackage{amsfonts}
\newtheorem{theorem}{Theorem}[section]

\newtheorem{proposition}[theorem]{Proposition}

\newcommand{\proof}{\noindent{\bf Proof.\ }}
\newcommand{\qed}{\hfill $\square$ \bigskip}
\def\cp{\,\square\,}

\begin{document}

\title{Some binary products and integer linear programming for computing $k$-metric dimension of graphs}

\author{Sandi Klav\v zar $^{a}$ \and Freydoon Rahbarnia $^{b,}$\footnote{Corresponding author} \and Mostafa Tavakoli $^b$} 

\date{}

\maketitle

\vspace*{-10mm}
\begin{center}

$^a$ Faculty of Mathematics and Physics, University of Ljubljana, Slovenia\\
{\tt sandi.klavzar@fmf.uni-lj.si}

\medskip



$^b$ Department of Applied Mathematics, Faculty of Mathematical Sciences,\\
Ferdowsi University of Mashhad, P.O.\ Box 1159, Mashhad 91775, Iran\\
{\tt rahbarnia@um.ac.ir}\\
{\tt tavakoli2002@gmail.com}

\end{center}

\begin{abstract}
Let $G$ be a connected graph. For an ordered set $S=\{v_1,\ldots, v_\ell\}\subseteq V(G)$, the vector $r_G(v|S) = (d_G(v_1,v), \ldots, d_G(v_\ell,v))$ is called the metric $S$-representation of $v$. If for any pair of different vertices $u,v\in V(G)$, the vectors $r(v|S)$ and $r(u|S)$ differ in at least $k$ positions, then $S$ is a $k$-metric generator for $G$. A smallest $k$-metric generator for $G$ is a $k$-{\em metric basis} for $G$, its cardinality being the $k$-metric dimension of $G$. A sharp upper bound and a closed formulae for the $k$-metric dimension of the hierarchical product of graphs is proved. 
Also, sharp lower bounds for the $k$-metric dimension of the splice and link products of graphs are presented.
An integer linear programming model for computing the $k$-metric dimension and a $k$-metric basis of a given graph is proposed. These results are applied to bound or to compute the $k$-metric dimension of some classes of graphs that are of interest in mathematical chemistry.
\end{abstract}


\medskip\noindent
{\bf AMS Subj.\ Class:} 05C12; 05C76

\section{Introduction} 
\label{sec:intro} 
Throughout this paper, all graphs are assumed to be connected, finite, and simple. Let $G = (V(G), E(G))$ be a graph and $u, v\in V(G)$. The distance $d_G(u,v)$ is defined as the length of a shortest path that connects $u$ and $v$. 
For an ordered set $S=\{v_1,\ldots, v_\ell\}\subseteq V(G)$, the vector $r_G(v|S) = (d_G(v_1,v), \ldots, d_G(v_\ell,v))$ is called the metric $S$-representation of $v$. The set $S$ {\em distinguishes} vertices $u$ and $v$ if $r_G(u|S)\neq r_G(v|S)$.
We say that $S$ is a {\em metric generator} for $G$ if for every pair of distinct vertices $u,v \in V(G), r_G(u|S)\neq r_G(v|S)$.
A metric generator with the smallest cardinality among all the metric generators for $G$ is
called a metric basis {\em metric basis} for $G$, and its cardinality the {\em metric dimension} of $G$, denoted by ${\rm dim}(G)$. 
 
The concept of the metric dimension was introduced almost half a century ago in~\cite{Harary, Slater}. Afterwards, the concept was exdtensively investigated, cf.~\cite{bailey-2011, Chartrand-2000, Das, Rahbarnia}, see also~\cite{johnson-1993, khuller-1996} for its applications in modeling of real world problems. Moreover, several versions of this concept such as the local metric dimension~\cite{okamoto-2010}, strong resolving sets~\cite{Oellermann-Peters-Fransen}, and edge metric dimension~\cite{kelenc-2018} were introduced, in particular because of their applications in modelling of different problems. 

As an extension of metric generators, $k$-metric generators were recently proposed in~\cite{Estrada-Moreno} as follows. If $G$ is a graph and $k$ a natural number, then $S\subseteq V(G)$ is a $k$-{\em metric generator} for $G$ if for any pair of different vertices $u,v\in V(G)$, there exist at least $k$ vertices $v_1, \ldots, v_k$ in $S$ such that $d_G(v_i,u)\neq d_G(v_i,v)$, for every $i\in [k]$, where $[k]$ denotes the set $\{1,\ldots,k\}$. We further say that if $X\subseteq V(G)$, then S is a $k$-{\em metric generator} for $X$ if any pair of vertices from $X$ is distinguished by at least $k$ vertices from $S$. Note that if $|X|\le 1$, then, by definition, $S = \emptyset$ is a $k$-{\em metric generator} for $X$. A smallest $k$-metric generator is a $k$-{\em metric basis} for $G$ and the $k$-{\em metric dimension} ${\rm dim}_k(G)$ of $G$ is the cardinality of a $k$-metric basis for $G$. If $G$ admits no  $k$-metric basis, then we set ${\rm dim}_k(G) = \infty$. Moreover, if $G$ admits a $k$-{\em metric generator}, then we will write ${\rm dim}_k(G) < \infty$.

Independently from~\cite{Estrada-Moreno}, the $k$-metric dimension was introduced and studied in~\cite{adar-2017}, where besides the unweighted version of the problem, also weighted version of the problem was studied for paths, complete graphs, complete bipartite graphs, and complete wheel graphs. (See also~\cite{kang-2019} for the fractional version of the $k$-metric dimension.) In~\cite{kmetric-NP-hardness} it is proved that the problem of computing the $k$-metric dimension of graphs is NP-hard and that the problem can be solved in linear time for some special trees. Moreover, it is proved that for a connected graph $G$, the problem of finding the largest integer $k$ such that  $G$ admists a $k$-metric generator can be solved in polynomial time. The $k$-metric dimension (more precisely $k$-resolving sets) was used in~\cite{bailey-2019} to construct  error-correcting codes. In the same paper the $k$-metric dimension of Cartesian products of paths was determined. The $k$-metric dimension of the lexicographic product of graphs was investigated in~\cite{kmetric-lexicographic} and the $k$-metric dimension of the corona product of graphs in~\cite{kmetric-corona}. Finally, in~\cite{corregidor-2021} related bounds were investigated, while in~\cite{beardon-2019} the concept was studied on general metric spaces. 

In this paper we proceed this line of investigation by considering the $k$-metric dimension of hierarchical, splice, and link products. In the next section we give a general sharp upper bound on ${\rm dim}_k(G(U)\sqcap H)$ and an exact formula for the case $|U| = 1$. 
Also, we give general sharp lower bounds on ${\rm dim}_k(G(U)\cdot H)$ and ${\rm dim}_k(G(U)\sim H)$. In Section~\ref{sec:ILP} we propose an integer linear programming model for computing the $k$-metric dimension of an arbitrary graph. In the last section we use results of Sections~\ref{sec:hierarchical-products} and~\ref{sec:ILP} to bound or determine the $k$-metric dimension of some classes of graphs that appear in mathematical chemistry.

\section{Some binary operations} 
\label{sec:hierarchical-products} 

In the first subsection, we study the $k$-metric dimension under the hierarchical product of graphs. In the second subsection, we investigate the $k$-metric dimension under the splice and the link product of graphs.

\subsection{Hierarchical products}

Suppose $G$ and $H$ are two graphs with $U\subseteq V(G)$. The {\em hierarchical product} $G$ and $H$ with respect to $U$, denoted by $G(U) \sqcap H$, is the graph with the vertex set $V(G) \times V(H)$, and the edge set $\{(g, h)(g', h')\ |\ g = g'\in U\  \text{and}\ hh' \in E(H), \ \text{or},\ gg'\in E(G)\ \text{and}\  h =h'\}$.
If $h\in V(H)$, then the subgraph of $G(U) \sqcap H$ induced by the vertices $(g,h)$, $g\in V(G)$, is isomorphic to $G$ and called a {\em $G$-layer}. Similarly {\em $H$-layers} are defined, cf.~\cite{klavzar}. We note that the operation $\sqcap$ (for two and also more factors) was in the seminal paper~\cite{barriere-2009} named the {\em generalized hierarchical product}. Also, we note that 
if $U=V(G)$, then $G(U) \sqcap H$ is the standard Cartesian product $G\cp H$, cf.~\cite{klavzar1}.
Moreover, if $|U|=1$, then $G(U) \sqcap H$ is a cluster product $G\{H\}$, see \cite{Schwenk}.
We refer to~\cite{anderson-2017, anderson-2018} for different studies of the hierarchical product. 

Suppose $u$ and $v$ are two vertices of $G$, and $U\subseteq V(G)$. A $u,v$-walk $W$ is a {\em $u,v$-walk through $U$} if $W$ is an $u,v$-walk in $G$ that contains some vertex of $U$. Note that the latter vertex could be one of $u$ and $v$. In the following, $d_{G(U)}(u,v)$ denotes the length of a shortest $u,v$-walk in $G$ through $U$. 

\begin{proposition} {\rm \cite{barriere-2009}}
\label{prp:distance} 
For two graphs $G$ and $H$ with $U\subseteq V(G)$, we have 
$$d_{G(U)\sqcap H}((g,h),(g',h'))=\begin{cases} 
d_{G(U)}(g,g')+d_H(h,h'); & h\neq h',\\ 
d_G(g,g'); & h=h'. 
\end{cases}$$ 
\end{proposition} 
Before stating our main results, we need to introduce some notations. For $v\in V(G)$ and $\ell\in {\mathbb N}_0$, we use the notation $N_\ell(v)$ to denote the set of vertices of $G$ that are at distance $\ell$ from $v$, that is, 
$$N_\ell(v) = \{u\in V(G):\ d_G(v,u) = \ell\}\,.$$
If $U\subseteq V(G)$, then we set: 
$${\rm dim}_k(G(U)) = \min\{ \mid\! \bigcup_{u\in U\atop \ell\ge 1} S_G(u,\ell)\! 
\mid\ :\ S_G(u,\ell)\ \text{is}\ k\text{-metric\ generator\ for}\ N_\ell(u) \}\,.$$
In words, ${\rm dim}_k(G(U))$ is the size of a smallest set of vertices $S$ such that each pair of vertices from each of the sets $N_\ell(u)$ is distinguished by $k$ vertices from $S$. If such a set does not exist, then we set  ${\rm dim}_k(G(U)) = \infty$.  Denoting the order of a graph $G$ by $n(G)$, we now have the following result. 

\begin{theorem} 
\label{t1} 
Let $G$ be a graph, $U\subseteq V(G)$ and ${\rm dim}_k(G(U))=t$. 
If $H$ is a graph with ${\rm dim}_{\lceil \frac{k}{t}\rceil}(H) < \infty$, then 
\[{\rm dim}_k(G(U)\sqcap H)\leq n(H)\cdot {\rm dim}_k(G(U)). \] 
\end{theorem} 

\proof 
We use $X$ instead of $G(U)\sqcap H$ for convenience. Let $S^T(G)\subseteq V(G)$ be a set that realizes ${\rm dim}_k(G(U))$, so that $|S^T(G)| = {\rm dim}_k(G(U))$. Let further $S_G(u,\ell)$ be a subset of $S^T(G)$ which is a $k$-metric generator for $N_u(\ell)$, where $u\in U$ and $\ell \ge 1$. 

Set $S=S^T(G)\times V(H)$.  
We will show that $S$ is a $k$-metric generator for $X$. To do this, let $v$ and $w$ be arbitrary, different vertices of $X$, and consider the following cases. 
 
\medskip\noindent 
{\bf Case 1}: $v=(g,h)$ and $w=(g,h')$, that is, $v$ and $w$ belong to a common $H$-layer.\\ 
Since ${\rm dim}_{\lceil \frac{k}{t}\rceil}(H) < \infty$, there exist $\lceil\frac{k}{t}\rceil$ vertices say $h_1, \ldots, h_{\lceil\frac{k}{t}\rceil}$ in $H$ such that $d_H(h,h_i)\neq d_H(h',h_i)$ for every $i\in [ \lceil\frac{k}{t}\rceil ]$. Set \[S'=\{(g',h_i)\; : \; g'\in S^T(G)\; \text{and}\; i\in [ \lceil\frac{k}{t}\rceil ] \}\, .\] 
Clearly, $|S'| = \lceil \frac{k}{t}\rceil\cdot t \geq k$, $S'\subseteq S$, and $d_X((g,h),(g',h_i))-d_X((g,h'),(g',h_i))=d_H(h,h_i)-(h',h_i)\neq 0$ for every $(g',h_i)\in S'$. 
 
\medskip\noindent 
{\bf Case 2}: $v=(g,h)$ and $w=(g',h)$, that is, $v$ and $w$ belong to a common $G$-layer.\\ 
If there exist $u\in U$ and $\ell\in {\mathbb N}$ such that $d_G(g,u)=d_G(g',u) = \ell$, then $S_G(u,\ell)\times\{h\}$ is a subset of $S \times \{h\}$ which is a $k$-metric generator for $\{v, w\}$. Otherwise, $d_G(g,u)\neq d_G(g',u)$ for each $u\in U$. Thus, for each $(g'',h'')\in S-(S^T(G)\times \{h\})$, we have $d_X((g,h),(g'',h'')) - d_X((g',h),(g'',h''))=d_{G(U)}(g,g'')-d_{G(U)}(g',g'')>0$. On the other hand, $|S-(S^T(G)\times \{h\})|\geq k$ because $|S^T(G)|\geq k$. 
 
\medskip\noindent 
{\bf Case 3}: $v=(g,h)$ and $w=(g',h')$, where $g\ne g'$ and $h\ne h'$.\\ 
If there exist $u\in U$ and $\ell\in {\mathbb N}$ such that $d_G(g,u)=d_G(g',u) = \ell$, then there exist $k$ vertices $g_1,\ldots,g_k\in S_G(u,\ell)$ with $d_G(g,g_i)\neq d_G(g',g_i)$ for every $i\in [k]$. Assume, without loss of generality, that $d_G(g,g_i)>d_G(g',g_i)$ for $i\in [q]$, and $d_G(g,g_i)<d_G(g',g_i)$ for $q<i\leq k$. Therefore, $d_X((g,h),(g_i,h'))=d_{G(U)}(g,g_i)+d_H(h,h')\geq d_G(g,g_i)+d_H(h,h')>d_G(g',g_i)=d_X((g',h'),(g_i,h'))$ for $i\in [q]$. Also, $d_X((g,h),(g_i,h))=d_G(g,g_i)<d_G(g',g_i)+d_H(h,h')\leq d_{G(U)}(g',g_i)+d_H(h,h')=d_X((g',h'),(g_i,h))$ for $q<i\leq k$. Then
$$\{(g_1,h'),\ldots,(g_q,h'),(g_{q+1},h),\ldots,(g_k,h)\}\,,$$ 
which is a subset of $S$, is a $k$-metric generator for $\{v, w\}$. 
 
Otherwise, $d_G(v,z)\neq d_G(w,z)$ for each $z\in U$. Consider a vertex $u\in U$ and an $\ell\in {\mathbb N}_0$. Suppose that $g_1,\ldots,g_k\in S_G(u,\ell)$. Assume, without loss of generality, that $d_G(g,g_i)\geq d_G(g',g_i)$ for $i\in [q]$, and $d_G(g,g_i)<d_G(g',g_i)$ for $q<i\leq k$. Therefore, $d_X((g,h),(g_i,h'))=d_{G(U)}(g,g_i)+d_H(h,h')\geq d_G(g,g_i)+d_H(h,h')>d_G(g',g_i)=d_X((g',h'),(g_i,h'))$ for $i\in [q]$. Also, $d_X((g,h),(g_i,h))=d_G(g,g_i)<d_G(g',g_i)+d_H(h,h')\leq d_{G(U)}(g',g_i)+d_H(h,h')=d_X((g',h'),(g_i,h))$ for $q<i\leq k$. Then 
$$\{(g_1,h'),\ldots,(g_q,h'), (g_{q+1},h),\ldots,(g_k,h)\}\,,$$ 
which is a subset of $S$, is a $k$-metric generator for $\{v, w\}$. 
 
We conclude that $S$ is a $k$-metric generator for $X$, and consequently  
${\rm dim}_k(X)\leq |S|=n(H)\cdot {\rm dim}_k(G(U))$. 
\qed 

In the last section, we will demonstrate sharpness of the bound of Theorem~\ref{t1}. Here we continue with the case $|U| = 1$. If $U = \{u\}$, then we simplify the notation $G(\{u\})$ to $G(u)$. If $G$ is a path and $u$ its end vertex, then we say $G(u)$ is a {\em rooted path}. 

\if 
Let $G(u)$ be a rooted graph. Then ${\rm edim}(G(u)\sqcap H)=1$ where $H$ be a path. Also, $\{(u,w)\ :\ w\in V_W\}$ is 
an edge metric generator for $G(u)\sqcap W$, and so ${\rm edim}(G(u)\sqcap W)\leq |V_W|$. 
Therefore, we can say that if $X = G(u)\sqcap H$, where $G(u)$ is a rooted path, then 
$1\leq {\rm edim}(X) \leq n(H)$. The next theorem gives the exact value of the edge metric dimension of hierarchical product of graphs $G(u)$ and $H$ where $G(u)$ is not a rooted path. 
\fi 

\begin{theorem} 
\label{t2} 
Suppose that ${\rm dim}_k(G(u))=t$, where $G(u)$ is not rooted path. If $H$ is a graph with ${\rm dim}_{\lceil \frac{k}{t}\rceil}(H) < \infty$, then
\[{\rm dim}_k(G(u)\sqcap H)= n(H)\cdot {\rm dim}_k(G(u)). \] 
\end{theorem} 

\proof 
Set again $X = G(u)\sqcap H$, and let $S^T(G)\subseteq V(G)$ be defined as in the proof of Theorem~\ref{t1}. In that theorem we have proved that $S=S^T(G)\times V(H)$ is a $k$-metric generator for $X$. To complete our proof we are going to show that $S$ is a $k$-metric basis of $X$. 

 Assume by contradiction that $S'$ is a $k$-metric generator for $X$ such that $|S'| < n(H)t$. Then, by the pigeonhole principle, there exists a vertex $h\in V(H)$ such that $|S'\cap (S^T(G)\times \{h\})|< |S^T(G)|$. This implies that there exist two vertices $(g',h),(g'',h)\in V(G)\times \{h\}$ such that $d_X((g',h),(u,h))=d_X((g'',h),(u,h))$. On the other hand, we have $d_X((g',h),w)=d_X((g',h),(u,h))+d_X((u,h),w)$ and $d_X((g'',h),w)=d_X((g'',h),(u,h))+d_X((u,h),w)$ for each $w\in S'-(V(G)\times \{h\})$. Thus, $d_X((g',h),w)=d_X((g'',h),w)$ for each $w\in S'-(V(G)\times \{h\})$. This means that $S'$ is not a $k$-metric generator for $X$, a contradiction. Therefore, each $k$-metric generator for $X$ has at least $n(H)\cdot {\rm dim}_k(G(u))$ vertices and $S$ is a $k$-metric basis of $X$. 
\qed 

\subsection{Splice and link products}

Let $G$ and $H$ be disjoint graphs and let $a\in V(G)$ and $b\in V(H)$.  The {\em splice} $(G\cdot H)(a; b)$ of $G$ and $H$ (with respect to $a$ and $b$) is the graph obtained from $G$ and $H$ by identifying the vertices $a$ and $b$. Similarly, the  {\em link} $(G\sim H)(a;b)$ of $G$ and $H$  (with respect to $a$ and $b$) is the graph obtained from $G$ and $H$ by adding the edge $ab$, cf.~\cite{SL1, SL2}. Simplifying the notation $G(\{a\})$ to $G(a)$ as already done above, we have then have the following lower bounds. 

\begin{theorem} \label{splic}
If $G$ and $H$ are disjoint graphs,  $a\in V(G)$, and $b\in V(H)$, then 
\begin{align*}
{\rm dim}_k((G\cdot H)(a;b)) & \geq {\rm dim}_k(G(a))+{\rm dim}_k(H(b))\,, \\
{\rm dim}_k((G\sim H)(a;b)) & \geq {\rm dim}_k(G(a))+{\rm dim}_k(H(b))\,.
\end{align*}
\end{theorem}

\proof
Set $X=(G\cdot H)(a;b)$. By $G'$ and $H'$, we denote the respective copies of $G$ and $H$ in $X$. Let $S$ be a $k$-metric generator for $X$ and assume to the contrary that  $|S| < {\rm dim}_k(G(a)) + {\rm dim}_k(H(b))$. Then $|S\cap V(G')|<{\rm dim}_k(G(a))$ or $|S\cap V(H')|<{\rm dim}_k(H(b))$. We may assume that $|S\cap V(G')|<{\rm dim}_k(G(a))$. Thus there exist vertices $u$ and $v$ in $G'$ such that  $d_X(u,a)=d_X(v,a)$ and they are not distinguished by at least $k$ vertices from $S\cap V(G')$. On the other hand,  since $d_X(u,a)=d_X(v,a)$, the vertices $u$ and $v$ cannot be distinguished by the vertices from $V(H')$, and in particular from  $S\cap V(H')$. Therefore, $S$ is not a $k$-metric generator for $X$. This contradiction proves the first inequality. 

The second inequality is proved along parallel lines. 
\qed

\begin{figure}[ht!] 
\centerline{ \includegraphics[scale=0.45]{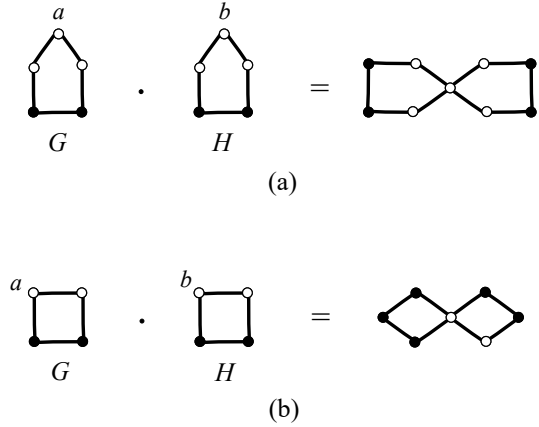}} 
\caption{Two splice products and their $2$-metric basis.} 
\label{fig4} 
\end{figure}

As an example, consider the $2$-metric dimension and the two examples as presented in Fig.~\ref{fig4}. The bold vertices in each of the graphs form its $2$-metric basis. The first example (Fig.~\ref{fig4}(a)) demonstrates that the bound of Theorem~\ref{splic} is sharp, while the second example (Fig.~\ref{fig4}(b)) shows that the equality need not hold. Indeed, one can verify that in the latter example we have ${\rm dim}_2(G(\{a\})) = {\rm dim}_2(H(\{b\})) = 2$ and ${\rm dim}_2((G\cdot H)(a;b))=5$.

One would be tempted to replace in Theorem~\ref{splic} the values ${\rm dim}_k(G(a))$ and ${\rm dim}_k(H(b))$ with ${\rm dim}_k(G)$ and ${\rm dim}_k(H)$, respectively. However, this cannot be done. For instance, since ${\rm dim}_2(K_n) = n$, $n\ge 3$, we infer that for $n,m\ge 3$,  
$${\rm dim}_2(K_n)+{\rm dim}_2(K_m) = n + m > n(K_n\cdot K_m)(a;b)) \geq {\rm dim}_2((K_n\cdot K_m)(a;b))\,.$$
Similarly, ${\rm dim}_2(P_n) = 2$ for $n\ge 2$. Setting $G=P_n$, $H=P_m$, and selecting $a$ and $b$ to be the pendant vertices of $G$ and $H$, respectively, we have $(G\cdot H)(a;b) = P_{n+m-1}$ and hence ${\rm dim}_2((G\cdot H)(a;b)) = 2$ as well. Note that this example also shows that ${\rm dim}_k(G(a))$ and ${\rm dim}_k(H(b))$ cannot be replaced with   ${\rm dim}_k(G)$ and ${\rm dim}_k(H)$ in the second inequality of Theorem~\ref{splic}. Note finally that if $a$ is a pendant vertex of $P_n$, then $ {\rm dim}_2(P_n(a)) = 0$ because there are no pairs of vertices in $P_n$ that are at the same distance from $a$.

\section{Integer LP model} 
\label{sec:ILP} 
In \cite{Chartrand-2000}, Chartrand et al.\ adopted an integer linear programming model (ILPM) to obtain the metric dimension and a metric basis for a graph. Motivated by this work, we present an ILPM for obtaining the $k$-metric basis as follows. 
Let $G$ be a graph with $V(G) = \{v_1,\ldots,v_n\}$. Let $D_G=[d_{ij}]$ be the distance matrix of $G$, that is, an $n\times n$ matrix with $d_{ij}=d_G(v_i,v_j)$, $i,j\in [n]$. For $x_i \in \{0,1\}$, $i\in [n]$, define the function $F$ by 
\[F(x_1,\ldots, x_n)=x_1+\cdots+x_n\,, \] 
and set
\[ \delta(d_{ij}, d_{i'j'}) = 
   \begin{cases}
          1; & d_{ij} \ne d_{i'j'}\,,  \\
          0; & d_{ij} = d_{i'j'}\,.
   \end{cases} 
\]
The the goal is to minimize $F$ subject to the constraints 
\[ \delta(d_{i1}, d_{j1}) x_1 + \cdots + \delta(d_{in}, d_{jn}) x_n \geq k,\ \ 1 \le i < j \le n\,.\] 
Note that if $x'_1, \ldots, x'_n$ is a set of binary values for which $F$ attains its minimum, then $W = \{v_i:\ x'_i =1\}$ is a $k$-metric basis for $G$. 

As a simple example consider the path $P_3$ on verices $v_1$, $v_2$, $v_3$. Then \[ D_{P_3}=\begin{pmatrix} 
0 & 1& 2\\ 
1 & 0 & 1\\ 
2 & 1 & 0 
\end{pmatrix}.\] 
Let $k=2$. Then the ILPM is to minimize $F(x_1,x_2,x_3)=x_1+x_2+x_3$ subject to the constraints $x_1+x_2+x_3 \geq 2$, $x_1 + x_3 \geq 2$, $x_1+x_2 +x_3 \geq 2$, $x_1, x_2, x_3 \in \{0, 1\}$. $F$ attains its minimum for $x_1=1$, $x_2=0$, and $x_3=1$, hence $W=\{v_1,v_3\}$ is a $2$-metric basis for $P_3$. 

\section{Applications} 
\label{sec:applications} 

In this section we first apply Theorem~\ref{t1} to obtain upper bounds on the $k$-metric dimension of some families of graphs that appear in chemical graph theory. Then, using the ILMP, we compute exact values for some smaller graphs from these families and conclude that the bound of Theorem~\ref{t1} is sharp in several examples. Throughout this section,  let $S^T(G)$ have the same meaning as in the proof of Theorem~\ref{t1}, that is, a set that realizes ${\rm dim}_k(G(U))$. 
 
\begin{figure}[ht!] 
\centerline{ \includegraphics[scale=0.35]{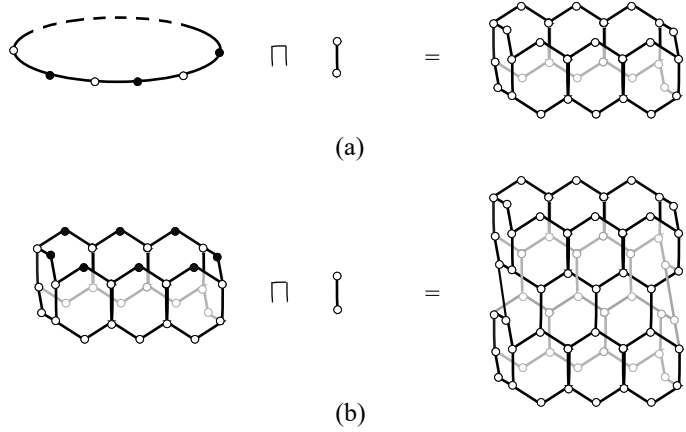}} 
\caption{ (a) $C_{2p}(U)\sqcap P_2=F_{p,1}$ where $U$ is formed by the bold vertices\ (b) $F_{p,1}(U)\sqcap P_2=F_{p,3}$ where $U$ is formed by the bold vertices.} 
\label{fig1} 
\end{figure}

Consider $F_{p,t}$, a {\em zigzag nanotube} with $t$ hexagonal belts with $p$ hexagons in each belt. 
As it is shown in Fig.~\ref{fig1}(a), $F_{p,1}$ is isomorphic to $C_{2p}(U)\sqcap P_2$ where $C_{2p}$ is a cycle with the ordered vertex set $\{v_1,\ldots,v_{2p}\}$ and $U=\{v_{2i}:\ i\in [p]\}$. Clearly, if $k\leq p$, then $S^T(C_{2p})=\{v_{2i-1}:\ i\in [k]\}$, and $S^T(C_{2p})=\{v_{2i-1}:\  i\in [p]\} \cup \{v_{2i}:\ i\in [k-p+1]\}$ otherwise. Thus, 
$${\rm dim}_k(C_{2p}(U)) = 
\begin{cases} 
k; & k\leq p\,,\\ 
k+1; & p<k<2p\,. 
\end{cases}$$  
 On the other hand, ${\rm dim}_1(P_2) < \infty$, hence by Theorem \ref{t1}, we have 
\[{\rm dim}_k(F_{p,1})={\rm dim}_k(C_{2p}(U)\sqcap P_2)\leq n(P_2)\cdot {\rm dim}_k(C_{2p}(U))=\begin{cases} 
2k; & k\leq p\,,\\ 
2k+2; & p<k<2p\,.
\end{cases}\] 
Now, consider $F_{p,3}$ depicted in Fig.~\ref{fig1}(b). The fact that ${\rm dim}_k(F_{p,1}(U))\leq {\rm dim}_k(F_{p,1})$ and Theorem~\ref{t1} lead us to the following bound for the $k$-metric dimension of $F_{p,3}$: 
\begin{align*} 
{\rm dim}_k(F_{p,3})&={\rm dim}_k(F_{p,1}(U)\sqcap P_2)\leq n(P_2)\cdot {\rm dim}_k(F_{p,1}(U))\\ &\leq n(P_2)\cdot {\rm dim}_k(F_{p,1})= 
\begin{cases} 
4k; & k\leq p\,,\\ 
4k+4; & p<k<2p\,. 
\end{cases}  
\end{align*} 
By induction we infer that 
\[{\rm dim}_k(F_{p,2q-1})\leq \begin{cases} 
2^qk; & k\leq p\,,\\ 
2^q(k+1); &  p<k<2p\,. 
\end{cases}\tag{1} \] 
 
\begin{figure}[ht!] 
\centerline{ \includegraphics[scale=0.31]{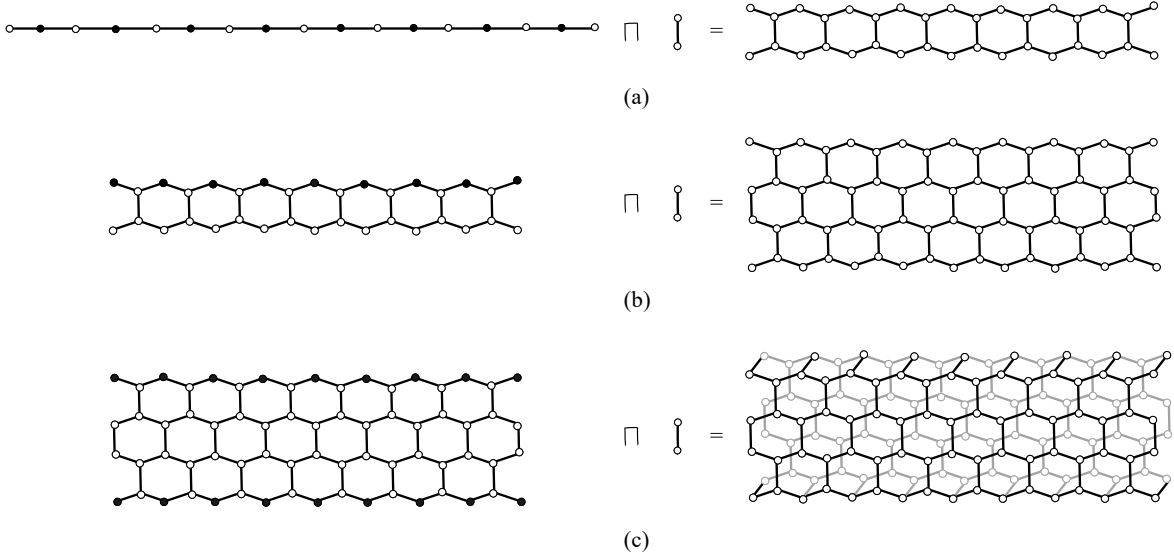}} 
\caption{A sequence of graphs $\Gamma_{1,7}$, $\Gamma_{3,7}$ and $A_{8,7}$ constructed by  the hierarchical product; the final graph $A_{8,7}$ is an armchair graph.} 
\label{fig2} 
\end{figure}

A {\em zig-zag polyhex lattice-like} $\Gamma_{2t-1,p}$ is a planar graph which is formed by $2t-1$ hexagonal rows with $p$ and $p+1$ hexagons in the rows, alternatively, and a pendent vertex at both ends of 
its first and last level. See $\Gamma_{1,7}$ and $\Gamma_{3,7}$ in Fig.~\ref{fig2}(a) and Fig.~\ref{fig2}(b), respectively. An {\em armchair graph} $A_{4t,p}$ is a tube whose surface is covered by $4t$ hexagonal rows with $p$ and $p+1$ hexagons in the rows, alternatively. The graph $A_{8,7}$ shown in Fig.~\ref{fig2}(c) is an armchair graph. 
 
Let $P_{2p+3}$ be a path with ordered vertex set $\{v_1,\ldots,v_{2p+3}\}$. The graph $\Gamma_{1,p}$ is isomorphic to $P_{2p+3}(U)\sqcap P_2$ where $U=\{v_{2i}:\ i\in [p]\}$ (Fig.~\ref{fig1}(a) shows $\Gamma_{1,7}=P_{2p+3}(U)\sqcap P_2$ where $U$ is formed by the black vertices in the figure). Easily one can check that if $k\leq p+2$, then $S^T(P_{2p+3})=\{v_{2i-1}:\ i\in [k]\}$, and otherwise
\[S^T(C_{2p})=\{v_{2i-1}:\  i\in [p+2]\} \cup \{v_{2i}:\ i\in [k-p-1]\}.\]
Thus, 
\[{\rm dim}_k(P_{2p+3}(U))=\begin{cases} 
k; & k\leq p+2\,,\\ 
k+1; & p+2<k<2p+3\,. 
\end{cases}\]
On the other hand, ${\rm dim}_1(P_2) < \infty$ and hence, applying Theorem~\ref{t1} again, 
\[ {\rm dim}_k(\Gamma_{1,p})\leq n(P_2)\cdot {\rm dim}_k(P_{2p+3}(U)) = 
\begin{cases} 
2k; & k\leq p+2\,,\\ 
2k+2; & p+2<k<2p+3\,. 
\end{cases} \tag{2} \] 
Then, by (2), \[{\rm dim}_k(\Gamma_{1,7})\leq 
\begin{cases} 
2k; & k\leq 9\,,\\ 
2k+2; & 9<k<17\,. 
\end{cases}\] 
Now consider $\Gamma_{3,7}$ depicted in Fig.~\ref{fig1}(b). By the fact that ${\rm dim}_k(\Gamma_{1,7}(U))\leq {\rm dim}_k(\Gamma_{1,7})$ and Theorem~\ref{t1}, we can estimate as follows: 
\begin{align*} 
{\rm dim}_k(\Gamma_{3,7})&={\rm dim}_k(\Gamma_{1,7}(U)\sqcap P_2)\leq n(P_2)\cdot {\rm dim}_k(\Gamma_{1,7}(U))\\ &\leq n(P_2)\cdot {\rm dim}_k(\Gamma_{1,7})= 
\begin{cases} 
4k; & k\leq 9\,, \\ 
4k+4; & 9<k<17\,.
\end{cases}  
\end{align*} 
For another example consider the armchair graph $A_{8,7}$ from Fig.~\ref{fig2}(c). This figure demonstrates that $A_{8,7}$ is isomorphic to $\Gamma_{3,7}(U)\sqcap P_2$, where $U$ is formed by the black vertices. Then, by Theorem~\ref{t1}, we have 
\begin{align*} 
{\rm dim}_k(A_{8,7})&={\rm dim}_k(\Gamma_{3,7}(U)\sqcap P_2)\leq n(P_2)\cdot {\rm dim}_k(\Gamma_{3,7}(U))\\ &\leq n(P_2)\cdot {\rm dim}_k(\Gamma_{3,7})= 
\begin{cases} 
8k; & k\leq 9\,, \\ 
8k+8; & 9<k<17\,. 
\end{cases}
\end{align*} 

\begin{table}[h!] 
\centering 
\renewcommand{\arraystretch}{2}{ 
\begin{tabular}{|c|c|c|c|} 
\hline 
  & $F_{4,1}$ & $\Gamma_{1,2}$ & $\Gamma_{1,3}$ \\ 
\hline\hline 
bound on $2$-metric dimension &$4$ & $4$ &$4$ \\ 
\hline 
exact value of $2$-metric dimension &$4$& $4$ &$4$ \\ 
\hline \hline 
bound on $3$-metric dimension &$6$ &$6$ & $6$ \\ 
\hline 
exact value of $3$-metric dimension &$6$ &$5$ & $5$  \\ 
\hline\hline 
bound on $4$-metric dimension &$8$ & $8$ & $8$ \\ 
\hline 
exact value on $4$-metric dimension &$8$ &$7$ &$7$\\ 
\hline\hline 
bound on $5$-metric dimension &$10$ & $10$ & $10$ \\ 
\hline 
exact value of $5$-metric dimension &$9$ &$8$ & $9$ \\ 
\hline 
\end{tabular}}
\caption{A comparison of the bounds obtained from relations (1) and (2) and the exact values of $k$-metric dimension obtained by ILPM .} 
\label{table:2} 
\end{table} 
 
A comparison of the bounds~(1) and~(2) and the exact values of $k$-metric dimension obtained by ILPM for the graphs $F_{4,1}$, $\Gamma_{1,2}$ and $\Gamma_{1,3}$, and $k\in \{2,3,4,5\}$ is shown in Table~\ref{table:2}. The information of this table shows that the bound presented in Theorem~\ref{t1} is sharp in several cases.

\begin{figure}[ht!] 
\centerline{ \includegraphics[scale=0.25]{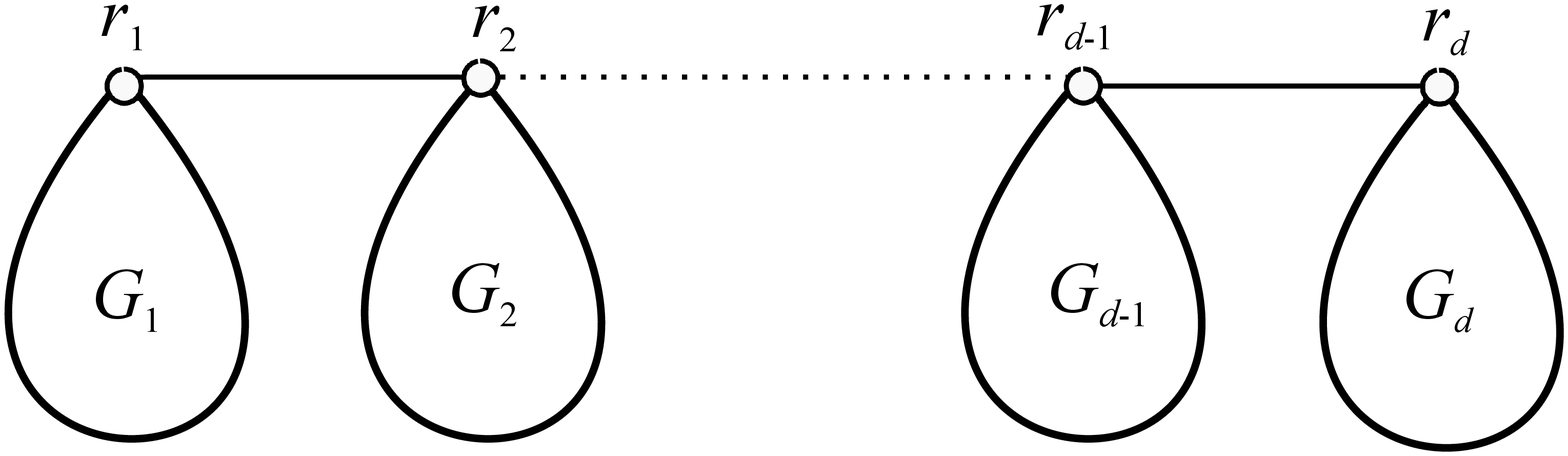}} 
\caption{The bridge-path graph $BP(G_1,\ldots,G_d; r_1,\ldots,r_d)$.} 
\label{fig3} 
\end{figure}

To conclude the paper we give an example in which Theorem~\ref{t2} is applied. Let $G_1,\ldots,G_d$ be graphs with $v_i\in V(G_i)$. The bridge-path graph of these graphs with respect to the vertices $r_1,\ldots,r_d$ is the graph $BP(G_1,\ldots,G_d; v_1,\ldots,v_d)$ obtained from the graphs $G_1,\ldots,G_d$ by connecting the vertices $r_i$ and $r_{i+1}$ by an edge for all $1\leq i\leq d-1$, see Fig.~\ref{fig3}. 

If $G_1=\cdots=G_d=G$ and $r_1=u\in V(G)$, then we have
\[BP(G_1,\ldots,G_d; r_1,\ldots,r_d)\cong G(u)\sqcap P_d.\]
Thus, if $G(u)$ is not a rooted path and ${\rm dim}_k(G(u))=t$, then by Theorem~\ref{t2} we have 
\[{\rm dim}_k(BP(G_1,\ldots,G_d; r_1,\ldots,r_d))={\rm dim}_k(G(u)\sqcap P_d)=d\cdot {\rm dim}_k(G(u)),\]
  where $\lceil\frac{k}{t}\rceil\leq d-1$.

\end{document}